\documentclass[12pt,reqno]{amsart}
\usepackage{}
\usepackage{bbm}
\usepackage{amsfonts}
\usepackage{mathrsfs}

\usepackage{amssymb}
\usepackage{amsmath}

\renewcommand{\phi}{\varphi}
\newcommand{\dd}{{\mathrm{d}}}

\newcommand{\R}{{\mathbb R}}

   \newcommand{\EX}{{\mathbb E}}
   \newcommand{\PX}{{\mathbb P}}

\newtheorem{theorem}{Theorem}[section]
\newtheorem{lemma}[theorem]{Lemma}

\newtheorem{corollary}[theorem]{Corollary}
 \numberwithin{equation}{section}

\theoremstyle{definition}
\newtheorem{definition}[theorem]{Definition}
\newtheorem{remark}[theorem]{Remark}

\begin{document}
\date{November 1, 2010 \\This work was partly supported by the   NSFC grants   10971225 and 11028102, the NSF Grant 1025422,  and the Cheung Kong Scholars Program. }
\title{Evolution systems of measures for stochastic flows }
\author{Xiaopeng ~Chen}
\address[X. Chen]{%
 School of Mathematics and Statistics\\ Huazhong University of Science and
Technology\\ Wuhan 430074, China} \email[X.~Chen]{chenxiao002214336@yahoo.cn}
\author{Jinqiao ~Duan}
\address[J. Duan]{%
Department of Applied Mathematics
\\ Illinois Institute of Technology\\ Chicago, IL 60616, USA \\
\& School of Mathematics and Statistics\\ Huazhong University of Science and
Technology\\ Wuhan 430074, China
} \email[J.~Duan]{duan@iit.edu}
\author{Michael ~Scheutzow}
\address[M. Scheutzow]{%
 Institut f\"ur Mathematik, MA 7-5, Technische Universit\"at Berlin, Stra{\ss}e des 17. Juni 136, D-10623 Berlin}
\email[M.~Scheutzow]{ms@math.tu-berlin.de}

\maketitle



\begin{abstract}

A new concept of {\em an evolution system of measures
for stochastic flows}  is considered. It corresponds to the notion of an invariant
measure for  random dynamical systems (or cocycles).  The existence of  evolution systems of
measures for asymptotically compact stochastic flows  is obtained.
For a white noise stochastic flow, there exists a one to one correspondence  between evolution systems of
measures for a  stochastic flow  \emph{and} evolution systems
of measures for the associated Markov transition semigroup. As an application,
an alternative  approach for evolution systems of measures of
 2D stochastic Navier-Stokes equations with  a time-periodic forcing term is presented.

\end{abstract}




\section{Introduction}

An {\em invariant measure} for a random dynamical system (RDS) defined on a probability space $(\Omega,\mathscr{F},\PX)$
taking values in a Polish space $X$ is a
probability measure $\mu$ on $\Omega\times
X$ that is invariant under the skew-product flow, and whose marginal
on $\Omega$ coincides with the given measure $\PX$ (see \cite{Arn,Cra}). Any
(forward) invariant random compact set supports at least one
invariant measure \cite{Arn,Cra}.  Moreover, there exists an invariant  measure for each asymptotically compact RDS \cite{Brz}.
Important examples of an RDS are those generated by stochastic differential equations (SDEs) (see \cite{Arn,AS}) and
certain stochastic partial differential equations (SPDEs).  On the
other hand, invariant measures for the Markov semigroup generated by an SDE or SPDE
have gained great attention,
see \cite{Pra} for a comprehensive survey. The two notions of an invariant measure are rather different.  The first is
from ergodic theory and dynamical systems, and the second is from probability theory.  A link between these two
notions is that of a Markov measure \cite{Cra3, Cra1}.

In applications to physics, an equation can be simultaneously subjected to  noise and time dependent non-stationary
forcing \cite{Pra3}.  In this case the theory of random dynamical systems does not apply and one is lead to consider
nonautonomous stochastic systems or {\em stochastic flows}.
The theory of nonautonomous  deterministic systems has been considered by many authors(e.g. \cite{Bon, Che, Har, Vis}).  Haraux \cite{Har} investigated
the concept of a {\em uniform attractor} paralleling that of a global attractor.
The theory of attractors for nonautonomous deterministic systems is considered in the book \cite{Vis}.

In this paper,  we introduce a new concept of {\em evolution systems of measures
for stochastic flows}, which -- in some sense -- generalize the notion of invariant
measures for random dynamical systems. We show that any attractor of a stochastic flow supports at least one such
evolution system
and in our main result, Theorem \ref{1.1}, we establish a relationship between evolution systems of measures
for {\em white noise} stochastic flows and evolution systems
of measures for transition operators which have been considered recently in the study of the long time behavior
of SPDEs \cite{Gei, Pra1,Pra2,Pra3}.
As an application,
we consider the stochastic Navier-Stokes equation subjected to a time dependent forcing
term and prove that it generates a white noise stochastic flow which admits an attractor and hence an evolution system of
measures (both for the flow and the Markov semigroup).


\section{Preliminaries}
We recall some definitions from \cite{Cra2}. Let
$(\Omega,\mathscr{F},\PX)$ be a
probability space. Further,
let $(X,d)$ be a Polish (i.e.~complete separable metric) space and denote the
Borel-$\sigma$-algebra of $X$ by $\mathscr{B}(X)$.
\begin{definition}\label{SF}
A family of maps
$S(t,s;\omega):X\rightarrow X$,
$-\infty <s \leq t<\infty$, $\omega \in \Omega$ is called a {\em stochastic flow (SF)}  if for $\PX$-a.e. $\omega$
\begin{itemize}
\item[(i)] $S(t,r;\omega)S(r,s;\omega)=S(t,s;\omega)$ for all $s\leq r\leq
t$;\\
\item[(ii)]  $x \mapsto S(t,s;\omega)x$ is continuous for all $s\leq t$;\\
\item[(iii)] For all $t\in \mathbb{R}$, $x\in X$ the mapping
$(s,\omega)\mapsto S(t,s;\omega)x$ is measurable from
$((-\infty,t]\times \Omega,\mathscr{B}((-\infty,t])\otimes
\mathscr{F})$ to $(X,\mathscr{B}(X))$;\\
\item[(iv)] $S(t,t;\omega)=\mathrm{id}_{X}$ for all $t \in \R$.
\end{itemize}
\end{definition}
Given $t\in \mathbb{R}$
and $\omega\in \Omega$,
we say that $K(t,\omega)\subset X$ is an
{\em attracting set} if, for every non-empty bounded set  $B\subset X$,
$$d(S(t,s;\omega)B,K(t,\omega))\rightarrow 0 \quad \mbox{as }\quad s\rightarrow
-\infty,$$
where $d(A,B):=\sup_{x \in A}\inf_{y \in B}d(x,y)$.
We say that the stochastic flow
$(S(t,s;\omega))_{t\geq s,\omega\in \Omega}$ is {\em asymptotically compact}
if there exists a measurable set $\Omega_0$ with measure one such
that for all $t\in \mathbb{R}$ and $\omega\in \Omega_0$, there
exists a compact attracting set $K(t,\omega)$ \cite{Cra2}.

We define the random omega limit set of a bounded set $B\subset X$ at time $t$ as
$$A(B,t,\omega)=\bigcap\limits_{T<t}\overline{\bigcup\limits_{s<T}S(t,s;\omega)B}$$
and define
$$A(t,\omega)=\overline{\bigcup\limits_{B\subset X}A(B,t,\omega)}.$$

For a measurable space $(\Omega,{\mathscr{F}})$,  we denote the family of all
probability measures on that space by ${\mathcal{P}}(\Omega)$. We denote by $\mathcal{P}_{\PX}(\Omega\times X)$
the set of all probability measures on
$(\Omega\times X,{\mathscr{F}}\otimes \mathscr{B}(X))$ with marginal $\PX$ on $\Omega$.

For $\mu\in \mathcal{P}_{\PX}(\Omega\times X)$, we denote its disintegration by
$\mu_\omega$, i.e.~for each $\omega \in \Omega$, $\mu_\omega \in \mathcal{P}(X)$, the map
$\omega \mapsto \mu_\omega (B)$ is measurable for all $B \in \mathscr{B}(X)$ and
$\mu(A \times B)=\int 1_A \mu_\omega(B) \dd \PX(\omega)$ for all $A \in {\mathscr{F}}$ and
$B \in {\mathscr{B}(X)}$.
Note that such a disintegration exists since the space $X$ is Polish.

Next we define an evolution system of measures for a stochastic flow.
\begin{definition}
An {\em evolution system of measures} for
a stochastic flow $S$ is a family of probability measures $\{\mu_s\}_{s \in \R}$ on
$(\Omega\times X,\, \mathscr{F}\otimes \mathscr{B}(X))$ which satisfies $$ S(t,s;.)\mu_s=\mu_t,$$ for each $s \leq t$
and whose marginals on $\Omega$ coincide with the given measure $\PX$.
\end{definition}
Note that $S(t,s;.)\mu_s$ stands for the image of the measure $\mu_s$ under the map $(\omega,x) \mapsto
(\omega,S(t,s;\omega)x)$.

For each $s\in \mathbb{R}$, denote $$\mathcal {F}_{\geq s}=\overline{\sigma\{S(u,t;\cdot)x\mid s\leq t\leq u
,x\in X\}}$$ and  $$\mathcal
{F}_{\leq s}=\overline {\sigma\{S(u,t;\cdot)x\mid t\leq u \leq s,x\in X\}},$$
where the bar denotes the completion with respect to the measure $\PX$.

\begin{definition}
We  call a stochastic flow $(S(t,s;\omega))_{t\geq
s,\omega\in \Omega}$  a  {\em white noise stochastic flow},
if for each $s\in \mathbb{R}$,  $\mathcal {F}_{\geq s}$ and $\mathcal
{F}_{\leq s}$ are independent.
\end{definition}
The standard example of a white noise SF is the flow generated by a stochastic differential equation
driven by Brownian motion.  If $(S(t,s;\omega))_{t\geq
s,\omega\in \Omega}$ is a  white noise SF, then we define $$P_{s t}f(x):=\mathbb{E}(f(S(t,s;\omega)x))$$ for $t\geq s$, and $f$
a bounded measurable function from $X$ to $\mathbb{R}$. Note that this family defines an (inhomogeneous) Markov
semigroup of transition operators, i.e. we have $P_{s u} =  P_{tu} \circ P_{s t}$ whenever $s \leq t \leq u$.

For $\rho \in \mathcal{P}(X)$, let
$P_{st}\rho$ denote the image of $\rho$ under $P_{st}$.
A family of Borel probability measures $\{\rho_s\}_{s \in \R}$ on $X$ with $P_{s t}\rho_s=\rho_t$ for all
$s\leq t$ is said to be an {\em evolution system of measures} for $\{P_{s t}\}$.

\section{Main result}

We quote the following result from reference \cite{Cra2} (Theorem 2.1).
\begin{lemma}\label{lemma}
  Assume  that  the  stochastic flow
$(S(t,s;\omega))_{t\geq s,\omega\in \Omega}$ is  asymptotically
compact and define $K$ and $A$ as in the previous section. Then,  for $\mathbb{P}$-a.e.~$\omega$,  the following  results hold
true. For  all  $t \in \mathbb{R}$, the  set $A( t, \omega)$  is a
nonempty  compact attracting subset  of $K(t, \omega)$,  and it  is  the minimal closed  set
with this property. Moreover, $A$  is  invariant,  in the  sense that
for all $s \leq t$, $$ S( t, s; \omega) A(s, \omega) =A(t,\omega).$$
In this case, $A$ is called  a  {\em pullback attractor} of the stochastic flow $S$.
\end{lemma}

We remark that under the conditions of the previous lemma, the set $A(t,\omega)$ is measurable with respect to the $\PX$-completion
of $\mathscr{F}$ for each $t \in \R$ (this is Proposition 2.1 in  \cite{Cra2}). An inspection of the proof shows that $A(t,\omega)$
is even measurable with respect to ${\mathcal{F}}_{\le t}$ for every $t \in \R$.\\


We now have the following result.

\begin{theorem}\label{1.2}
Let $S$ be a stochastic flow and suppose that $A(t,\omega)$, $t \in \R$, $\omega \in \Omega$  is a compact forward invariant
family of sets, i.e.~there exists a set of full measure $\Omega_0 \in \mathscr{F}$ such that
\begin{itemize}
\item[a)] for each $x \in X$, $t \in \R$ the map $\omega \mapsto d(x,A(t,\omega))$ is Borel measurable,
\item[b)] for each $\omega \in \Omega_0$,  $t \in \R$ the set $A(t,\omega)$ is compact,
\item[c)] for each $\omega \in \Omega_0$,  $s \le t$ we have $S(t,s;\omega)A(s,\omega)\subseteq A(t,\omega)$.
\end{itemize}
Then there exists an evolution system of measures $\{\mu_t\}_{t\in \R}$ for $(S(t,s;\omega))_{t\geq s,\omega\in \Omega}$
such that $\mu_t$ is supported by $A(t,\omega)$, i.e. for each $t \in \R$ we have $\mu_{t \omega}(A(t,\omega))=1$ for $\PX$-a.e. $\omega$.
\end{theorem}
\begin{proof}
We  assume that $\Omega_0$ is chosen such that properties (i), (ii) and (iv) of Definition \ref{SF} hold for all
$\omega \in \Omega_0$. For $\omega \in \Omega_0$ we define
$$
C_{n,m}(\omega):=S(-m,-n;\omega)A(-n,\omega),\quad m = 0,1,2,...,\quad n=m,m+1,...
$$
Then
\begin{align*}
C_{n+1,m}(\omega)&=S(-m,-(n+1);\omega)A(-(n+1),\omega)\\
&=S(-m,-n;\omega) S(-n,-(n+1);\omega)A(-(n+1),\omega)\\
&\subseteq S(-m,-n;\omega) A(-n,\omega)= C_{n,m}(\omega), \qquad n = m,m+1,...
\end{align*}
Since all sets $C_{n,m}(\omega)$ are nonempty and compact the same is true for
$$
\bar C_m(\omega):=\cap_{n=m}^{\infty} C_{n,m}(\omega), \qquad m=0,1,2,...
$$
and $\bar C_m$ enjoys the measurability property a) in the theorem, i.e.~$\bar C_m$ is a random compact set.
The measurable selection theorem of Castaing and Valadier \cite{CV}, Theorem III.9, guarantees the existence of a measurable
map $x_0$ from $\Omega_0$ to $X$ such that $x_0(\omega) \in \bar C_0(\omega)$ for all $\omega \in \Omega_0$. By induction,
we find a measurable map $x_{-(n+1)}$ from $\Omega_0$ to $X$ such that $x_{-(n+1)}(\omega) \in \bar C_{n+1}(\omega)$ and
$S(-n,-(n+1);\omega) x_{-(n+1)}(\omega) = x_{-n}(\omega) $ for all $\omega \in \Omega_0$ (again using the measurable selection theorem).
Finally, for $s \in \R$
we define $x_s(\omega):=S(s,-n;\omega)x_{-n}(\omega)$ for some nonnegative integer $n$ such that $-n \le s$
(by the previous construction the
definition does not depend on the choice of $n$). Then we have $x_t(\omega)=S(t,s;\omega)x_s(\omega)$ for all
$\omega \in \Omega_0$ and all $-\infty<s\le t < \infty$ which implies that $\mu_{t\omega}:=\delta_{x_t(\omega)}$ is an
evolution system supported by $A$.
\end{proof}
From Lemma \ref{lemma} and Theorem \ref{1.2} we have the following result.
\begin{corollary}
There existence of  a  evolution systems of
measure  for asymptotically compact stochastic flows.
\end{corollary}
\begin{remark}
For an RDS which has a pullback attractor, it is well-known that
the attractor supports all invariant measures \cite{Cr}. It is generally not true however that an attractor of an SF
supports all evolution systems of measures. As an example consider the flow which is generated by
the ODE $x'(t)=-x(t)$, $t \in \R$ on $X=\R$ which has an attractor $A(t,\omega)=\{0\}$ which does not support the
evolution system of measures $\mu_t:=\delta_{\exp\{-t\}},\,t \in \R$.\\
\end{remark}

\begin{remark}
There exist stochastic flows which do not admit an evolution system of measures $\{\mu_{t}\}$. As a (deterministic) example take
a probability space consisting of a singleton,
$X=[0,\infty)$ equipped with the Euclidean metric and $S(t,s)x:=x+t-s,\;x \in X,\;-\infty<s \le t <\infty$.
\end{remark}


We will need the following elementary lemma.
\begin{lemma}\label{1.5}
  Let $(\Omega,\mathscr{F},\mathbb{P})$ be a probability space and $\mathscr{G}$ a sub $\sigma$-algebra of $\mathscr{F}$.
  Let $(E,d)$ be a
  Polish space with Borel $\sigma$-algebra $\mathscr{E}$. Suppose $f:\Omega\times E\rightarrow \mathbb{R}$ is a bounded
  $\mathscr{G} \otimes \mathscr{E}$-measurable function and $\mu_\omega$ a random
  probability measure on $E$ which is independent of $\mathscr{G}$, i.e.\\
  (i) $\omega\mapsto \mu_{\omega}(A)$ is $\mathscr{F}$-measurable and independent of $\mathscr{G}$ for each $A\in \mathscr{E}$.\\
  (ii) $\mu_{\omega}(\cdot)$ is a probability measure for each $\omega\in \Omega$. \\
  Define the probability measure $\rho$ on $(E,\mathscr{E})$ by
  $$\rho(A):=\mathbb{E}\mu_{\omega}(A), \quad\quad A\in \mathscr{E}. $$
  Then, almost surely
 \begin{eqnarray}\label{1.4}
 \mathbb{E}(\int f(\omega,x)\mu_{\omega}(\dd x)\mid \mathscr{G})=\int f(\omega,x)\rho(\dd x).
 \end{eqnarray}
 \end{lemma}
  \begin{proof}
    The statement holds for $f(\omega,x):=1_A(x)1_B(\omega)$, $A\in \mathscr{E} $, $B\in \mathscr{G}$.
    Both sides of \eqref{1.4} equal $1_B(\omega)\rho(A)$, so for $U\in \mathscr{G} \otimes \mathscr{E}$,
    $$ \mathbb{E}(\int 1_U(\omega,x)\mu_{\omega}(\dd x)\mid \mathscr{G})=\int 1_U(\omega,x)\rho(\dd x).$$
    Then the conclusion follows by linearity and the monotone convergence theorem.
  \end{proof}


Now we are ready to prove our main result. We wish to point out that the construction of an evolution system for the flow
from the evolution system for the Markov semigroup (part (a)) is largely analogous to the construction of an invariant measure
for a white noise RDS given an invariant measure for the Markov semigroup, see \cite{Arn}, Section 1.7 and \cite{Cra4}.
We provide a complete proof since the one in \cite{Arn} for RDS is rather sketchy and ours seems to be easier than
the one in \cite{Cra4}. In addition, there are some subtleties to be observed in our set-up due to the lack of continuity assumptions
with respect to the temporal variable which require us to restrict the convergence in \eqref{Gleichung} to deterministic sequences.

\begin{theorem}\label{1.1}
Suppose $(S(t,s;\omega))_{t\geq s,\omega\in \Omega}$ is a white noise
stochastic flow with the associated Markov semigroup $\{P_{st}\}$. Then
\begin{itemize}
\item[(a)] For an  evolution system of measures $\{\rho_s\}$ for $\{P_{s t}\}$ and any $t \in \R$,
\begin{equation}\label{Gleichung}
\lim\limits_{s\rightarrow-\infty}S(t,s;\omega)\rho_s=\mu_{t\omega}
\end{equation}
exists in the sense of weak convergence almost surely along each deterministic sequence $\{ s\} $ converging to $-\infty$
and

\begin{itemize}
\item[(i)] $\{\mu_t\}$ is an evolution system of measures for $S$. \\
\item[(ii)] $\mu_{t \omega}$ is measurable with respect to
            $\mathcal{F}_{\leq t}$ for each $t \in \R$.\\
\item[(iii)] $\EX (\mu_{t\omega}(B))=\rho_t(B)$ for all $t\in \mathbb{R}$, $B\in \mathscr{B}(X)$.\\
\end{itemize}
\item[(b)] If $\{\mu_t\}$ is an evolution system of  measures for $S$ such that  $\mu_t$ is $\mathcal
{F}_{\leq t}$-measurable for each $t\in \mathbb{R}$, then
$$\rho_t(B)=\int _\Omega \mu_{t \omega}(B)\dd \PX(\omega), \qquad t \in \R,\;B \in \mathscr{B}(X)$$
defines an evolution system of measures for $\{P_{st}\}$.
\end{itemize}
\end{theorem}
\begin{proof} (a) Fix $t\in \mathbb{R}$ and a bounded and measurable function $f$. We claim that the process
$$M_{s}(\omega):=\int f(x)(S(t,t-s;\omega)\rho_{t-s})(\dd x)$$ is a martingale with
respect to the filtration $\mathcal {G}_{s}:=
\mathscr{F}_{\geq t-s}$, $s \ge 0$. Let  $0\leq s \leq u$. Then by Lemma \ref{1.5} we have
\begin{eqnarray*}
  \mathbb{E}(M_{u}\mid \mathcal
{G}_{s})&=&\mathbb{E}\Big(\int f(x)(S(t,t-u;\omega)\rho_{t-u})(\dd x)\mid
\mathcal
{G}_{s}\Big)\\&=&\mathbb{E}\Big(\int f(S(t,t-u;\omega)y)\rho_{t-u}(\dd y)\mid
\mathcal
{G}_{s}\Big)\\&=&\mathbb{E}\Big(\int f(S(t,t-s;\omega)\circ S(t-s,t-u;\omega) y)\rho_{t-u}(\dd y)\mid
\mathcal
{G}_{s}\Big)\\&=&\mathbb{E}\Big(\int f(S(t,t-s;\omega)z)(S(t-s,t-u;\omega) \rho_{t-u})(\dd z)\mid
\mathcal
{G}_{s}\Big)\\&=&\int f(S(t,t-s;\omega)z) \rho_{t-s})(\dd z)\\
&=&\int f(x)(S(t,t-s;\omega)\rho_{t-s})(\dd x)=M_s,\quad a.s.
\end{eqnarray*}
Since $(M_s)$ is a bounded martingale, it converges both in $L^1$ and almost surely along any deterministic subsequence
to a random variable $X_f$ (since we did not impose continuity assumptions with respect to the temporal variables there is
no guarantee that $\lim_{s} M_s$ exists almost surely). In particular, $\EX X_f = \int f \dd \rho_t$.
Fix a sequence $s_n \to \infty$ and let $\nu_n(\omega):= S(t,t-s_n;\omega)\rho_{t-s_n}$. It is straighforward to check that
the sequence $\{\nu_n\}$ is almost surely tight. If $\nu$ is an almost sure limit point of the sequence, then
$\int f \dd \nu = X_f$ almost surely for each fixed bounded and continuous $f$ implying that $\nu$ is unique almost surely,
so the sequence $\{\nu_n\}$ converges almost surely to $\nu$ (in the topology of weak convergence of measures). Since $\nu$ does not
depend on the choice of the sequence $(s_n)$, \eqref{Gleichung}  follows (by identifying $\nu$ and $\mu_t$).
Further, we have $\EX \int f \dd \nu=\int f \dd \rho_t$ for each bounded and continuous $f$, so $\EX \mu_t=\rho_t$, i.e. (iii) follows.
It remains to verify (i) and (ii). \\

(i) Fix $-\infty< s\leq t<\infty$. Then  $\mu_{s\omega}=\lim_{u\rightarrow -\infty}S(s,u;\omega)\rho_{u}$,
$\mu_{t\omega}=\lim_{u\rightarrow -\infty}S(t,u;\omega)\rho_{u}$ for almost all $\omega$. For $u \leq s$ we have
$S(t,s;\omega)S(s,u;\omega)=S(t,u;\omega)$. Invoking the continuity of $S(t,s;\omega)$, we get
$$\mu_{t\omega}=\lim_{u\rightarrow -\infty}S(t,s;\omega)S(s,u;\omega)\rho_{s}=S(t,s;\omega)\mu_{s\omega}$$
almost surely, so (i) follows.\\

(ii)  $\mu_{t\omega}$ is $\mathscr{F}_{\leq t}$-measurable by construction. \\


(b) Let $\{\mu_s\}$ be an adapted evolution system of  measures for the flow $S$.  For every Borel set $B$ and $s \leq t$,
we have
\begin{align*}
 \rho_t(B)=\mathbb{E}\mu_{t \omega}(B)&=\mathbb{E}((S(t,s;\omega)\mu_{s\omega})(B))\\
&=\int\int1_B(S(t,s;\omega)(x))\mu_{s\omega}(\dd x)\PX (\dd \omega).
\end{align*}
If we can show that
\begin{eqnarray}\label{1.3}
\int\int f(\omega,x)\mu_{s\omega}(\dd x)\PX (\dd \omega) =\int\int f(\omega,x)\rho_{s}(\dd x) \PX (\dd \omega)
\end{eqnarray}
for each $\mathcal {F}_{\geq s}\otimes \mathscr{B}(X)$-measurable $f$, then we have $\rho_t(B)=(P_{s t}\rho_s)(B)$. For
$f(\omega,x):=1_{F\times A}$,
where $F\in \mathcal
{F}_{\geq s}$ and $B\in \mathscr{B}(X)$, we have
\begin{eqnarray*}
\mathbb{E}\int 1_{F\times A}\mu_{s \omega}(\dd x)=\mathbb{E}(1_{F}\mu_{s\omega}(A))=\PX (F)\rho_s(A)=
\mathbb{E}\int 1_{F\times A}\rho_{s }(\dd x)\PX(\dd \omega),
\end{eqnarray*}
so the measures $\mu_{s\omega}(\dd x)\dd \PX$ and $\dd \rho_{s}\dd \PX$ both restricted to $\mathcal
{F}_{\geq s}\otimes \mathscr{B}(X)$ coincide and so \eqref{1.3} follows.
\end{proof}


\begin{remark}
Theorem \ref{1.1} does not define a one-to-one correspondence between evolution systems for $\{S(t,s;\omega)\}$ and $\{P_{st}\}$
in general.
As an example take the identity flow $S$ on any Polish space $X$ which contains (at least) two distinct points $x_1$, $x_2$ defined on
a probability space containing a set $A \in \mathscr{F}$ such that $\PX(A)=1/2$. Then for each $\alpha \in [0,1]$
$$
\mu_{t\omega}:=
\begin{cases} \alpha \delta_{x_1} + (1-\alpha) \delta_{x_2}, \qquad \mbox{ if } \omega \in A\\
\\(1-\alpha) \delta_{x_1}+ \alpha \delta_{x_2}, \qquad \mbox{ if } \omega \notin A
\end{cases}
$$
is an evolution system of measures for the flow  and
$$
\rho_t(B):=\mathbb{E}\mu_{t\omega}(B)=\frac 12 \big(1_B(x_1)+1_B(x_2)\big)
$$
for $B\in  \mathscr{B}(X) $.   $\rho_t$, $t \in \R$ is an evolution system of measures (in fact even an invariant measure)
of the Markov semigroup of transition operators generated by $S$.
\end{remark}


 \section{An application}


In this section, we consider the 2D stochastic Navier-Stokes equations  subjected to a time-dependent (deterministic) forcing term as well as white   noise.
The existence of an evolution system of measures corresponding to the associated Markov semigroup has been investigated recently
\cite{Pra3}. The Theorems in the last section provide an alternative approach to this result. Namely,
we first show that the time-dependent 2D Navier-Stokes equations with additive noise generates
a white noise stochastic flow which is asymptotically compact and therefore, by Lemma \ref{lemma}, has a random
attractor. Hence   by Theorem \ref{1.2}, there exists an evolution system of measures for the stochastic flow which, by Theorem \ref{1.1},
generates an evolution system of measures for the associated Markov semigroup.

 More specifically, we consider the 2D Navier-Stokes equation  as follows:
 \[ \begin{cases} \dd u(t,\xi)=[\nu\bigtriangleup u(t,\xi)-(u(t,\xi)\cdot\nabla)u(t,\xi)]\dd t
 -\nabla p(t,\xi)\dd t+f(t,\xi)\dd t\\\qquad\quad\qquad+ \sum_{j=1}^m\phi_j\dd w_j(t), \\
\mbox{div}\, u(t,\xi)=0, \\  u(t,\xi)=0, \quad \xi\in \partial D,
\end{cases}
\]
where $D\subset \mathbb{R}^2$ is a bounded open domain with regular boundary  \cite{Cra2, Tem}, $w_j(t)$'s are independent two-sided real-valued Wiener processes on a probability space $(\Omega,\mathscr{F},\PX) $. Here $u$ is velocity, $\nu$ is viscosity, and  $p$ is pressure. The time-dependent forcing $f$ and noise intensities $\phi_j$'s are specified below.

Define $$H=\{u\in (L^2(D))^2\mid \mbox{div}\, u=0, u\cdot n=0 \ \mbox{on}  \ \partial D \} .$$
 The space $H$ is endowed with the usual scalar product $(\cdot,\cdot)$ and the associated norm $|\cdot|$.
Another useful space is
 $$V=(H_0^1(D))^2\bigcap H.$$
The space  $V$
 is endowed with the scalar product
 $$<u,v>=\sum_{i,j=1}^2(\frac {\partial u_i}{ \partial x_j}, \frac  {\partial v_i} { \partial x_j})$$  and the norm $$\|u\|^2=\sum\limits_{i,j=1,2}|\frac {\partial {u_i}}{\partial {x_j}}|_{L^2(D)}. $$
 The orthogonal projector in $(L^2(D))^2$ on the space $H$ is denoted by $\mathscr{P}$ and we define the Stokes operator
 $$A=-\mathscr{P} \Delta$$
 on
 $$D(A)=(H^2(D))^2\bigcap V.$$
 The bilinear operator $B:V\times V\rightarrow V'$ is given by $$\langle B(u,v),z\rangle= \int z(\xi)(u(\xi)\cdot\nabla)v(\xi)d\xi, \quad \quad u,v,z\in V.$$ See also \cite{FlaGat, Fla1} for more background information.

We assume that the noise intensities $\phi_j\in D(A)$, $1\leq k \leq m$, and there exists a constant $\beta>0$ such that
\begin{eqnarray*}
|<B(u,\phi_k), u>|\leq \beta |u|^2\quad \quad \mbox{ for all } u\in H, \ k=1,2,\cdots, m.
\end{eqnarray*}
The time-dependent forcing  $f:\mathbb{R}\longrightarrow
H$ is taken to be continuous and $2\pi$-periodic in time.

 We rewrite the Navier-Stokes equations in $H$ as:
  \begin{eqnarray*} \dd u(t)+\nu Au(t)\dd t+B(u(t),u(t))\dd t=f(t)\dd t+\sum_{j=1}^m \phi_j \dd w_j.
\end{eqnarray*}
Now use the change of variable
$$v(t)=u(t)-z(t),$$
where $z(t)=\sum_{j=1}^m\phi_jz_j$ is the Ornstein-Uhlenbeck process with
$$z_j=\int_{-\infty}^t e^{-\alpha(t-s)}dw_j(s),$$
with   a positive constant $\alpha$ to be   determined below.
It is known that $z(t)$ is a stationary  ergodic process and its trajectories are $\PX$-a.s.~continuous.

The new ``velocity" function $v(t)$ satisfies the following evolutionary equation with random coefficients

\begin{eqnarray}
  \frac {\dd v(t)} {\dd t} +\nu A v(t)+B(v(t)+z(t),v(t)+z(t))=f(t)+\alpha z-\nu A z\label{1.6}
  \end{eqnarray}
  and the initial condition
  \begin{eqnarray}
  v_s = u_s-z_s.  \label{1.7}
\end{eqnarray}

For each $\omega\in \Omega$,  by the Galerkin method,  it is known that for all $s\in \mathbb{R}$ and $v_s\in H$ with $v_s = u_s-z_s$, there exists a unique weak solution
$v(t,\omega)$ such that the mapping
$v_s \mapsto v(t, \omega; s,v_s)$ is continuous for all $t \geq s $ and $v(t,\cdot, s,v_s) $ is measurable.

We thus define the stochastic flow $(S(t,s;\omega))_{t\geq s,\omega\in \Omega}$ by
$$S(t,s;\omega)u_s=v(t,\omega)+z(t,\omega). $$
We now check that $S$ is indeed a  stochastic flow. \\
(i) \begin{eqnarray*}
  S(t,r;\omega)S(r,s;\omega)u_s&= &S(t,r;\omega)(v(r,\omega)+z(r,\omega))\\
  &= &  S(t,r;\omega)u_r=v(t,\omega)+z(t,\omega)=S(t,s;\omega)u_s.
\end{eqnarray*}
(ii)   $S(t,s;\omega)$ is continuous by the continuous $v(t,\omega)$ and $z(t,\omega)$.\\
 (iii) $S(\cdot,s;\cdot)u_s$ is measurable from  $((-\infty,t]\times \Omega,\mathscr{B}((-\infty,t])\times
\mathscr{F})$ to $(X,\mathscr{B}(X))$.

We now prove the existence of a compact attracting set $K(t, \omega)$ at time $t$.  The proof is similar to that given by Crauel, Flandoli and Debussche \cite{Cra,Cra2},  and we  only briefly describe some key ideas.  Let $B$ be
 a bounded set in $H$ and let $v$ be the solution of \eqref{1.6}-\eqref{1.7}. Multiplying  Equation \eqref{1.6}
  by $v$, we obtain

\begin{eqnarray*}
\frac 1 2  \frac {\dd |v|^2} {\dd t} +\nu\|v\|^2+(B(v+z,v+z),v)=(f,v)+\alpha(z,v)-\nu<Az,v>.\label{1.8}
  \end{eqnarray*}
  Note that
  \begin{eqnarray*}
    |(B(v+z,v+z),v)|&=&|(B(v+z,z),v+z)|\\
    &\leq& 2\beta(\sum_{j=1}^m|z_j|)(|v|^2+|z|^2).
  \end{eqnarray*}
It is also known that
  $$|u|\leq \lambda_1^{-1/2}\|u\|,  \quad \forall \ u\in V,$$
  where $\lambda_1$ is the first eigenvalue of $A$.
Thus it follows that
  \begin{eqnarray}\label{4.3}
  \frac {\dd |v|^2} {\dd t} +\frac \nu 4 \|v\|^2+(\frac {\nu \lambda_1} 4-2\beta\sum_{j=1}^m|z_j|)|v|^2\leq 2 g,\label{1.9}
  \end{eqnarray}
  where $g(t)$ is a function depending on $f$ and $z$.

  By the Gronwall inequality and the fact that $z$ is stationary and ergodic, for $s_1\leq s_0(\omega)$, $t_1\in [-1+t,t]$,
    \begin{eqnarray*}\label{a}
  |v(t_1)|^2\leq 2\beta_1 |u(s_1)|^2\exp (s_1 \frac {\nu\lambda_1} 8)+2 \beta_1 |z(s)|^2 \exp(s_1 \frac {\nu\lambda_1} 8)+ \nonumber\\
  2\beta_1\int_{-\infty}^t g(\sigma)\exp(\sigma( \frac {\nu\lambda_1} 4+\frac {2\beta} { \sigma} \int_{\sigma}^t\sum_{j=1}^m |z_j(\tau)|d\tau))d\sigma.
  \end{eqnarray*}
  So there exists $T(\omega,B)$,  depending only on   $B$ and $\omega$ such that for $s<T(\omega,B)$ and $t_1\in [-1+t,t]$,
    \begin{eqnarray*}v(t_1)\leq  2\beta_1\int_{-\infty}^t g(\sigma)\exp(\sigma( \frac {\nu\lambda_1} 4+\frac {2\beta} { \sigma} \int_{\sigma}^t\sum_{j=1}^m| z_j(\tau)|d\tau))d\sigma\\
    +2\beta_1\sup_{s\in(-\infty,-1]}(|z(s)|^2\exp(s \frac {\nu\lambda_1} 8))+1:= r_0(\omega).  \end{eqnarray*}
  Integrating \eqref{4.3} on $[-1+t,t]$ leads to
  \begin{eqnarray}
  \int_{-1+t}^t\| v(s)\|^2ds &\leq &\frac {8\beta }{\nu}( \int_{-1+t}^t\sum_{j=1}^m| z_j(\sigma)|d\sigma)r_0(\omega)+\frac 8 \nu \int_{-1+t}^t g(\sigma)d\sigma \nonumber \\
&  := &r_1(\omega).
  \end{eqnarray}
 Now     multiplying $v $ to Equation \eqref{1.6} in $V$,

\begin{eqnarray*}\label{2.1}
 \frac 1 2  \frac {\dd \|v\|^2} {\dd t} +\nu|Av|^2=<f,v>-(B(v+z,v+z),Av)+\alpha<z,v>-\nu(Az,Av).
  \end{eqnarray*}
 Therefore
  \begin{eqnarray*}
  \frac {\dd \|v\|^2} {\dd t}\leqslant G(t)+H(t)\|v\|^2,
  \end{eqnarray*}
  where $G(t)$ and $H(t)$ depend on $f$ and $z$.
  We deduce that for any $t_1\in [-1+t,t]$,
     \begin{eqnarray*}
\|v(t)\|^2\leq (\|v(t_1)\|^2+\int_{-1+t}^t G(\sigma)d\sigma )\exp(\int _{-1+t}^tH(\sigma)d\sigma).
  \end{eqnarray*}
  After integration on $[-1+t,t]$,
       \begin{eqnarray*}
\|v(t)\|^2\leq (\int_{-1+t}^t\|v(s)\|^2ds+\int_{-1+t}^t G(\sigma)d\sigma )\exp(\int _{-1+t}^tH(\sigma)d\sigma).
  \end{eqnarray*}
 Hence there exists $r(\omega)$ such that
     \begin{eqnarray*}
\|v(t)\|^2\leq r(\omega).
  \end{eqnarray*}
  Let $K$ be the ball in $V$ of radius $r(\omega)^{\frac 1 2}+\|z(t,\omega)\|$. Then  for any bounded set $B$ in $H$,
there exists $T(\omega,B)$ such that for $s<T(\omega,B)$,
  $$S(t,s;\omega)B\subset K(t,\omega).$$
 Then from the reference \cite[Theorem 2.1]{Cra2} there exists a random attractor. From Theorem \ref{1.2} we know that there exists an evolution system of measures.
Since the stochastic flow is a white noise stochastic flow,  there exists  a  corresponding  evolution system of measures for the Markov transition semigroup   by Theorem \ref{1.1}.

	


\bigskip
\bigskip

\end{document}